# Models of Prime-Like Sequences Generated by Least Element Sieve Operations Like the Sieve of Eratosthenes
## Friday, April 14, 2017

Leonard E. Baum
stefi_baum_80@post.harvard.edu

Abstract

We suggest other models of sieve generated sequences like the Sieve of Eratosthenes to explain randomness properties of the prime numbers, like the twin prime conjecture, the lim sup conjecture, the Riemann conjecture, and the prime number theorem.

Paper

The prime numbers 2,3,5,7,11,13,17,19,23,29,31,37… are the atoms of multiplication of the positive integers. They are indecomposable and every other positive integer is a product of prime numbers in a unique way. The prime numbers have numerous interesting asymptotic properties both global and local. For example, one of the global properties that has been proved is the so called prime number theorem which says that $\pi(n)$ defined as the number of prime numbers less than or equal to $n$ is asymptotic to *li(n)* defined as

$$li(n) \equiv \int_2^n 1/\ln(t)\, dt$$

. Another global property is the so-called Riemann conjecture which says that the real part of all the non-trivial complex zeros of $\zeta(s)$ equals 1/2 where

$$\zeta(s) = \sum_{n=1,\infty} 1/n^s = \prod_{primes\ p} (1 - \frac{1}{p^s})^{-1}$$

. The beautiful Riemann conjecture, which is equivalent to

$$\forall \epsilon > 0, \quad |\pi(n) - li(n)| = o(n^{(\frac{1}{2}+\epsilon)})$$

, is not really about analytic functions (and possibly is not provable by analytic methods). Hopefully, a stronger Brownian-motion type result like

$$\limsup(|\pi(n) - li(n)|/\sqrt{2\sum_{i=1,\ldots,n}\left[\frac{1}{\ln(i)\left(1-\frac{1}{\ln(i)}\right)}\right]}) = 1$$

is true. I suggest that the truth of the Riemann conjecture is due to the "numerical random" properties of the prime numbers as are certain local properties conjectured about the prime numbers such as the following two conjectures. One is the twin prime conjecture that there exists arbitrarily large n such that
$p_{n+1} - p_n = 2$. Another is that:
$\limsup(p_{n+1} - p_n) = (\ln p_n)^2$

These conjectures are due to certain randomness properties of the prime numbers that arise because the prime numbers are obtainable by the Sieve of Eratosthenes. I define the Sieve of Eratosthenes in the following specific way. At stage 0 we have the integers 1,2,3,4,5... At stage 1 we start the sieving process, and sieve the integers by the operator $\widehat{x_1} = 2$ removing multiples of 2 and leaving the odd numbers. At stage 2 we choose the smallest remaining number of the previous stage not equal 1 called $x_2 = 3$ and sieve by the operator $\widehat{x_2}$, removing all multiples of 3 thus leaving 1,5,7,11,13,17,19... At stage 3 we choose the smallest remaining integer not equal 1, which is $x_3 = 5$, and sieve by knocking out all multiples of 5 of the previous stage (removing 25, 35,55, 65, 85,95, ...). We then choose the smallest remaining integer not equal 1 and proceed in similar fashion knocking out multiples of 7, and continue this process infinitely many times. All integers $x_n$ that are un-sieved in any previous round (and then used as operators for sieving) are primes and those are all the primes. It is because the primes are obtainable by the sieving method that they have random qualities.

I now give a slightly altered model of the Sieve of Eratosthenes that will give an alternate set of integers as putative prime like objects. We start as before with integers 1,2,3,4,5.... At stage 1 we start just like before sieving by $x_1 = 2$, taking out every second element. At stage 2 we choose the smallest remaining not equal 1 integer $x_2 = 3$ and sieve out every third element being left with 1,5,7,11,13,17,19,23,25... as in the previous case. At the next stage we choose the smallest remaining not equal 1 unsieved element $x_3 = 5$, and sieve out every fifth remaining element leaving 7,11,13, and 17 but removing 19 and continue by leaving in 23,25,29, and 31 and taking out 35 etc. We then choose the smallest remaining not equal 1 element $x_4 = 7$, and leave in the next six elements, sieving out the 7th element and leaving in 11,13,17,23,25,29 and sieving out 31 proceeding through all the remaining integers. When sieving with 7 is complete we choose the smallest remaining not equal 1 element ($x_5 = 11$,) and sieve by that. And so on. In this model the "$x_n$ primes" are 1,2,3,5,7,11,13,17,23,25,29,37... and do not include the natural primes of 19 and 31 and do include the composite number 25. Do these $x_n$ primes have the same properties as the normal primes, such as: (i) the global prime number theorem (see Equation #1 below) and (ii) the Riemann conjecture (see Equation #2 below,), and (iii) the local generalized twin prime property in the "0-1" form that any finite constellation $a_1, \ldots, a_r$ of positive integer separations between $x_n\ x_{n+1}, \ldots, x_{n+r}$ occurs for infinitely many $n$ (with correct decreasing probability) unless it fails to occur at all in some finite stage, and (iv) $\limsup(x_{n+1} - x_n) = (\ln x_n)^2$ ?

Results for the real primes which rely on the primes being sufficiently plentiful and sufficiently random, like the Dirichlet theorem on the infinitude of primes in arithmetic sequences, or the Tao Green theorem that there are arbitrarily large arithmetic sequences all of which members are prime, or the Vinigardov theorem that all sufficiently large odd integers are sums of three odd not equal to one primes or the Goldbach conjecture that all even integers are sums of two odd not equal to one primes, should also hold for primes of the model.

The Sieve of Eratosthenes and the model yield examples of infinite sequences generated by rules (the divisibility and counting) of sieving out integers by least element operators $\widehat{x_n}$ in which each operator eliminates $1/x_n$ of the sequence on which it operates, so that at the end of the *n-th* stage a sequence with density,
d_n= $\prod_{i=1,n}(1 - 1/x_i)$ remains.  Observe that the model operators $\widehat{x_n}$ (perhaps the simplest of the class of such operators) will eliminate $1/x_n$ of any increasing infinite sequence on which it operates, while the divisibility rule of the Sieve of Eratosthenes operators $\widehat{x_n}$ require a special sequence of the result of the $\widehat{x_{n-1}}$ operator to succeed.

In the Sieve of Eratosthenes, and in the model, we keep in the integer 1, since no operator sieves it out. Also we use the hard sieve in which the sieving element is sieved out along with all the infinitely many elements it eliminates. Both procedures are necessary for the integers remaining at the end of each stage, *n*, to be periodic, with period $\prod_{i=1,n} x_i$. This periodicity follows by induction for the model. For the real primes of the Sieve of Erathosthenes, for i<=n, $x_i$ divides $\prod_{j=1,n} x_j$ + t if and only if $x_i$ divides t.

Conjecture: rules which leave periodic sequences with period equal to $\prod_{i=1,n} x_i$ and density $d_n$ equal to $\prod_{i=1,n}(1 - 1/x_i)$ ) after each stage n, obey the prime number theorem $x_n \to n\ln(n)$ (Rules may differ in smaller order terms such as $n \ln\ln(n)$ Periodicity should not be needed. It should be sufficient that the density of the infinite sequence left after stage n converges to $d_n$ sufficiently rapidly for each *n*.

For the model if we define $L_n(t)$ to be the number of integers <=t remaining after sieving by $\widehat{x_n}$ is complete at the end of stage *n*, since 1+[($L_{(n-1)}$ (t) – 1)/$x_n$)] are removed in stage n, we have

$L_n(t)$ = -1+(1-1/$x_n$) $L_{(n-1)}$+{ ($L_{(n-1)}$ (t)-1)/ $x_n$} and hence

$L_n(t)$ = -n + $\prod_{i=1,n}(1 - 1/x_i)t$ + $E_n$

where, $E_n = \prod_{i=1,n}(d_n/d_i)$ { $L_{i-1}(t)/ x_i$} with [] and {} denoting integer and fractional part, respectively.

Since for the model, $L_n(x_{n+i})$ = 1+i for $0 \le i \le x_n$ - 1,
the proof of Briggs and Hawkins (ref 4) for the lucky numbers of Ulam (ref 5), proving a Mertens theorem $d_n \sim 1/\ln(n)$) and a prime number theorem $x_n \sim n\ln(n)$ extends mutatis mutandis to prove Mertens theorem and the prime number theorem for the model. Unfortunately, this short, completely elementary method, does not extend to the sequence of real primes of Eratosthenes because the number of elements removed during stage n is not simply described by the least element $x_n$ and the number of elements <=t left after stage *n-1*.

For the model the relevant prime number theorem can be written as

**Equation #1**  $(x_{n+1} - 2 - \sum_{j=1,n} 1/d_j) / (2 + \sum_{j=1,n} 1/d_j) = o(1)$

The relevant Riemann like conjecture is that

**Equation #2**  $(x_{n+1} - 2 - \sum_{j=1,n} 1/d_j) = o\left(n^{\left(\frac{1}{2}+\epsilon\right)}\right) for\ all\ \epsilon > 0$

The stronger Brownian motion type behavior should also hold. However, Equations #1 and #2 are not true for the real primes of Erathosthenes because the expected waiting time for the least element

$x_{n+1}$ is not $1/d_n$

In fact by the real Mertens theorem and the real prime number theorem

$(1/n) \sum_{i=1,\ldots,n}[(x_{i+1} - x_i)d_i - 1] \rightarrow (e^{-\gamma}-1) \sim -0.44$

For the real primes of the Sieve of Eratosthenes with the divisibility rule if $x_n=m$ using the Mertens' theorem and the Chinese remainder theorem, Erdos (7) and Rankin (8) specify an integer *z* via residues modulo the primes $x_1 \ldots x_n$ such that each of the next

$k \cdot m \cdot \ln m \cdot \ln \ln \ln m / (\ln \ln m)^2$

integers following *z* is divisible by at least one of the primes $x_1 \ldots x_n$ and hence has been sieved out when the nth stage of sieving is complete. If *p* is the first prime following this interval, we see that the gap between *p* and the preceding prime is greater than

$k \cdot \ln p \cdot \ln \ln p \cdot \ln \ln \ln \ln p / (\ln \ln \ln p)^2$ .

This stage n result is much weaker than the conjectured Cramer result (ref 9) limsup $(p_{n+1} - p_n) = \ln^2(p_n)$.

There is a pseudo Chinese remainder theorem for the model. At the end of stage n, every remaining element is the $a_{i,j}$·th surviving element after an $x_i$ sieved out element. For $1 \leq a_{i,j} \leq x_i - 1$, for each i = 1, …, n and all such sequences, { $a_{i,j}$ } occur exactly once for survivors $\leq \prod_{i=1,n}(x_i)$. An induction argument also proves that the survivors after each stage n form a periodic sequence with period $\prod_{i=1,n}(x_i)$, as in the case of the Sieve of Erathosthenes. However, while the Mertens' theorem for the model allows elimination of large numbers of model survivors, as in the Erdos proof of gaps for the real primes, our pseudo Chinese remainder theorem for the model does not have good additive properties, so we do not have have Erdos Rankin like result for gaps after stage n for the model.

As a third case consider the Hawkins Sieve (ref 3) where the operator $\hat{x}_n$ in the n-th stage sieves out elements independently with probability $1/x_n$. David Williams and coauthors (ref 1) (ref 2) showed that the Strassen Theorem (ref 6) implies that with probability 1 the sums $\sum_{i=1,\ldots,n}[(x_{i+1} - x_i)d_i - 1]$ are quite close to Brownian motion. Then by a double use of Abel summation they conclude that, with probability 1, the error L $li\left(\frac{x_n}{L}\right) - n$ in the Prime Number Theorem (L a sequence dependent constant) is (omitting a relatively negligible amount) quite Brownian, a result stronger than the relevant Riemann conjecture for the Hawkins Sieve:

$$\left| L \; li\left(\frac{x_n}{L}\right) - n \right| = o(n^{\frac{1}{2}+\varepsilon}) \text{ for all } \varepsilon > 0 \text{ a.e.}$$

In the probabilistic model, where randomness is introduced externally in the sieving of each operator $\widehat{x_n}$, it is not surprising that a strong randomness arises for the least element sequence $x_1, x_2, x_3, \ldots$, and that this can be proved by good random variable techniques. In the deterministic case of the Sieve of Eratosthenes or the model when the operator $\widehat{x_n}$ eliminates every $x_n$-th element remaining after the $n-1$ stage, while one might guess the least element sequence $x_1, x_2, x_3, \ldots$ is numerically random, this seems more difficult to prove.

For the real primes, the enhanced Reiman conjecture is a statement about the random behavior (Brownian motion behavior) of the error π(n) - $li(n)$ in the estimate of the number of primes.

For any rule, the infinite sequence remaining after sieving by least element operators $\widehat{x_1}, \ldots, \widehat{x_n}$ is not random, being completely described by the rule and finitely many numbers $x_1, \ldots, x_n$. However, one may hope that the diagonal like infinite sequence of numbers $x_1, x_2, x_3, \ldots$ is numerically random.

Which rules for sieving the previous stage with the smallest integer remaining in that stage will obey the theorems and conjectures of the real primes and the model due to the numerical randomness of the sieve generated least elements?  Some regularity in the elements sieved out by each operator $\widehat{x_n}$ would be helpful.  Cases without regularity which are likely not to work are two part rules where $\widehat{x_n}$ sieves out everything (nothing) for increasingly very large finite times depending on $n$, and for the infinite remaining times sieve with density $1/x_n$. If in the model for $n \geq 3$ we alter the operation of $\widehat{x_n}$ so that it sieves out all elements until it comes to an integer $e$ such that $e+2$ is not next and then proceeds through the remaining integers by eliminating every $x_n$th element, then clearly the least element sequence $x_n$ will contain no twin primes for $n \geq 3$ although each sequence remaining after stage n will presumably contain infinitely many twin primes. (The lim sup$(x_{n+1} - x_n)$ situation is also altered.) Nevertheless, in this situation the least element sequence should still be quite "randomish." A good rule is non-anticipatory in the sense that like the sieve of Eratosthenes or the model, the decision to sieve out any element or not depends only on the information in that element or earlier elements of the sequence being acted on. The challenges are to give a suitable definition of "numerical random" so that "numerical random sequences" will have among others the four asymptotic local and global properties defined above and to define n-stage sieving operations so that the sequence of n-stage least elements will be numerically random in that sense.

Here are two additional questions to consider in the context of the standard primes of Eratosthenes. Because the Euler Riemann Zeta Function can be written as the Euler product over primes, properties of primes can be translated as properties of the Zeta function, and visa versa.  Question 1. What property of the primes is equivalent to the functional equation of the Zeta function?  Question 2. In the other direction, what property of the Zeta function would be equivalent to a Brownian motion-like behavior of  $\pi(n)$ - *li(n)* ?